\documentclass[11pt]{article}

\pdfoutput=1
\usepackage[pdftex]{color}
\usepackage{amssymb}
\usepackage{amsthm}
\usepackage{amsmath}
\usepackage{latexsym}
\usepackage{amscd}
\usepackage{graphicx}
\usepackage[pdftex, colorlinks=true, citecolor=green]{hyperref}
\usepackage{lscape}
\usepackage{setspace}
\usepackage{multirow}

\newcommand{\ds}{\displaystyle}

\newcommand{\ben}{\begin{equation}}     
\newcommand{\eeqn}{\end{equation}}
\newcommand{\bey}{\begin{eqnarray}}
\newcommand{\eey}{\end{eqnarray}}

\newtheorem{thm}{Theorem}[section]
\newtheorem{prop}[thm]{Proposition}

\newtheorem{defn}[thm]{Definition}
\newtheorem{conj}[thm]{Conjecture}

\begin{document}

\begin{flushleft}
{\Large
\textbf{Real and complex behavior for networks of coupled logistic maps}
}
\\
\vspace{4mm}
 Anca R\v{a}dulescu$^{*,}\footnote{Assistant Professor, Department of Mathematics, State University of New York at New Paltz; New York, USA; Phone: (845) 257-3532; Email: radulesa@newpaltz.edu}$, Ariel Pignatelli$^2$
\\
\indent $^1$ Department of Mathematics, SUNY New Paltz, NY 12561

\indent $^2$ Department of Mechanical Engineering, SUNY New Paltz, NY 12561
\\
\end{flushleft}

\vspace{3mm}
\begin{abstract}
\noindent Many natural systems are organized as networks, in which the nodes interact in a time-dependent fashion. The object of our study is to relate connectivity to the temporal behavior of a network in which the nodes are (real or complex) logistic maps, coupled according to a connectivity scheme that obeys certain constrains, but also incorporates random aspects. We investigate in particular the relationship between the system architecture and possible dynamics. In the current paper we focus on establishing the framework, terminology and pertinent questions for low-dimensional networks. A subsequent paper will further address the relationship between hardwiring and dynamics in high-dimensional networks.\\

\noindent For networks of both complex and real node-maps, we define extensions of the Julia and Mandelbrot sets traditionally defined in the context of single map iterations. For three different model networks, we use a combination of analytical and numerical tools to illustrate how the system behavior (measured via topological properties of the \emph{Julia sets}) changes when perturbing the underlying adjacency graph. We differentiate between the effects on dynamics of different perturbations that  directly modulate network connectivity: increasing/decreasing edge weights, and altering edge configuration by adding, deleting or moving edges. We discuss the implications of considering a rigorous extension of Fatou-Julia theory known to apply for iterations of single maps, to iterations of ensembles of maps coupled as nodes in a network. 
\end{abstract}

\section{Introduction}
\subsection{Network architecture as a system parameter}

\noindent Because many natural systems are organized as networks, in which the nodes (be they cells, individuals, populations or web servers) interact in a time-dependent fashion -- the study of networks has been an important focus in recent research. One of the particular points of interest has been the question of how the hardwired \emph{structure} of a network (its underlying graph) affects its \emph{function}, for example in the context of optimal information storage or transmission between nodes along time. It has been hypothesized that there are two key conditions for optimal function in such networks: a well-balanced adjacency matrix (the underlying graph should appropriately combine robust features and random edges) and well-balanced connection strengths, driving optimal dynamics in the system. However, only recently has mathematics started to study rigorously (through a combined graph theoretical and dynamic approach) the effects of configuration patterns on the efficiency of network function, by applying graph theoretical measures of segregation (clustering coefficient, motifs, modularity, rich clubs), integration (path length, efficiency) and influence (node degree, centrality). Various studies have been investigating the sensitivity of a system's temporal behavior to removing/adding nodes or edges at different places in the network structure, and have tried to relate these patterns to applications to natural networks.

Brain functioning is one of the most intensely studied contexts which requires our understanding of the tight inter-connections between system architecture and dynamics. The brain is organized as a ``dynamic network,'' self-interacting in a time-dependent fashion at multiple spacial and temporal scales, to deliver an optimal range for biological functioning.  The way in which these modules are wired together in large networks that control complex cognition and behavior is one of the great scientific  challenges of the 21st century, currently being addressed by large-scale research collaborations, such as the Human Connectome Project. Graph theoretical studies of empirical empirical data support certain generic topological properties of brain architecture, such as modularity, small-worldness, the existence of hubs and ``rich clubs''~\cite{bullmore2009complex,sporns2011non,sporns2002graph}.

In order to explain how connectivity patterns may affect the system's dynamics (e.g., in the context of stability and synchronization in networks of coupled neural populations), and thus the observed behavior, a lot of effort has been thus invested towards formal modeling approaches, using a combination of  analytical and numerical methods from nonlinear dynamics and graph theory, in both biophysical models~\cite{gray2009stability} and simplified systems~\cite{siri2007effects}. These analyses revealed a rich range of potential dynamic regimes and transitions~\cite{brunel2000dynamics}, shown to depend as much on the coupling parameters of the network as on the arrangement of the excitatory and inhibitory connections~\cite{gray2009stability}. The construction of a realistic, data-compatible computational model has been subsequently found to present many difficulties that transcend the existing methods from nonlinear dynamics, and may in fact require: (1) new analysis and book-keeping methods and (2) a new framework that would naturally encompass the rich phenomena intrinsic to these systems -- both of which aspects are central to our proposed work.

In a paper with Dr. Verduzco-Flores~\cite{radulescu2015nonlinear}, one of the authors of this paper first explored the idea of having network connectivity as a bifurcation parameter for the ensemble dynamics in a continuous time system of coupled differential equations. We used configuration dependent phase spaces and our probabilistic extension of bifurcation diagrams in the parameter space to investigate the relationship between classes of system architectures and classes of their possible dynamics, and we observed the robustness of the coupled dynamics to certain changes in the network architecture and its vulnerability to others. As expected, when translating connectivity patterns to network dynamics, the main difficulties were raised by the combination of graph complexity and the system's intractable dynamic richness. 

In order to break down and better understand this dependence, we started to investigate it in simpler theoretical models, where one may more easily identify and pair specific structural patterns to their effects on dynamics. The logistic family is historically perhaps the most-studied family of maps in nonlinear dynamics, whose behavior is by now relatively well understood. Therefore, we started by looking in particular at how dynamic behavior depends on connectivity in networks with simple logistic nodes. This paper focuses on definitions, concepts and observations in low-dimensional networks. Future work will address large networks, and different classes of maps.

Dynamic networks with discrete nodes and the dependence of their behavior on connectivity parameters have been previously described in several contexts over the past two decades. For example, in an early paper,  Wang considered a simple neural network of only two excitatory/inhibitory neurons, and analyzed it as a parameterized family of two-dimensional maps, proving existence of period-doubling to chaos and strange attractors in the network~\cite{wang1991period}. Masolle, Attay et al. have found that,  in networks of delay-coupled logistic maps, synchronization regimes and formation of anti-phase clusters depend on coupling strength ~\cite{masoller2011complex} and on the edge topology (characterized by the spectrum of the graph Laplacian)~\cite{atay2004delays}. Yu has constructed and studied a network wherein the undirected edges symbolize the nodes' relation of adjacency in an integer sequence obtained from the logistic mapping and the top integral function~\cite{yu2013logistic}.

In our present work, we focus on investigating, in the context of networked maps, extensions of the Julia and Mandelbrot sets traditionally defined for single map iterations. For three different model networks, we use a combination of analytical and numerical tools to illustrate how the system behavior (measured via topological properties of the \emph{Julia sets}) changes when perturbing the underlying adjacency graph. We differentiate between the effects on dynamics of different perturbations that  directly modulate network connectivity: increasing/decreasing edge weights, and altering edge configuration by adding, deleting or moving edges. We discuss the implications of considering a rigorous extension of Fatou-Julia theory known to apply for iterations of single maps, to iterations of ensembles of maps coupled as nodes in a network.

\subsection{Networking logistic maps}

The logistic map is historically perhaps the best-known family of maps in nonlinear dynamics. Iterations of one single quadratic function have been studied starting in the early 19th century, with the work of Fatou and Julia. 

The prisoner set of a map $f$ is defined as the set of all points in the complex dynamic plane, whose orbits are bounded. The escape set of a complex map is the set of all points whose orbits are unbounded. The Julia set of $f$ is defined as their common boundary $J(f)$. The filled Julia set is the union of prisoner points with their boundary $J(f)$.

For polynomial maps, it has been shown that the connectivity of a map's Julia set is tightly related to the structure of its critical orbits (i.e., the orbits of the map's critical points). Due to extensive work spanning almost one century, from Julia~\cite{julia1918memoire} and Fatou~\cite{fatou1919equations} until recent developments~\cite{branner1992iteration,qiu2009proof}, we now have the following: \\

\noindent {\bf Fatou-Julia Theorem.} {\it For a polynomial with at least one critical orbit unbounded, the Julia set is totally disconnected if and only if all the bounded critical orbits are aperiodic.}\\

\noindent For a single iterated logistic map\cite{carleson1993complex,devaney2006criterion}, the Fatou-Julia Theorem implies that the Julia set is either totally connected, for values of $c$ in the Mandelbrot set  (i.e., if the orbit of the critical point 0 is bounded), or totally disconnected, for values of $c$ outside of the Mandelbrot set (i.e., if the orbit of the critical point 0 is unbounded). In previous work, the authors showed that this dichotomy breaks in the case of random iterations of two maps~\cite{pignatelli2016}. In our current work, we focus on extensions for networked logistic maps. Although Julia and Mandelbrot sets have been studied somewhat in connection with coupled systems~\cite{isaeva2010phenomena}, none of the existing work seems to address the basic problems of how these sets can be defined for networks of maps, how different aspects of the network hardwiring affect the topology of these sets and whether there is any Fatou-Julia type result in this context. 

These are some of the questions addressed in this paper, which is organized as follows: In Section~\ref{logistic_maps}, we introduce definitions of our network setup, as well as of the extensions of Mandelbrot and Julia sets that we will be studying. In order to illustrate some basic ideas and concepts, we concentrate on three examples of 3-dimensional networks, which differ from each other in edge distribution, and whose connectivity strengths are allow to vary. In Section~\ref{complex_maps}, we focus on the behavior of these 3-dimensional models when we consider the nodes as complex iterated variables. We analyze the similarities and differences between node-wise behavior in each case, and we investigate the topological perturbations in one-dimensional complex slices of the Mandelbrot and Julia sets, as the connectivity changes from one model to the next, through intermediate stages. In Section~\ref{real_maps}, we address the same questions for real logistic nodes, with the advantage of being able to visualize the entire network Mandelbrot and Julia sets, as 3-dimensional real objects. In both sections, we conjecture weaker versions of the Fatou-Julia theorem, connecting points in the Mandelbrot set with connectivity properties of the corresponding Julia sets. Finally, in Section~\ref{discussion}, we interpret our results both mathematically and in the larger context of network sciences. We also briefly preview future work on high-dimensional networks and on networks with adaptable nodes and edges.


\section{Our models of networked logistic maps}
\label{logistic_maps}

\noindent We consider a set of $n$ nodes coupled according to the edges of an oriented graph, with adjacency matrix $\ds A=(A_{jk})_{j,k =1}^n$ (on which one may impose additional structural conditions, related to edge density or distribution). In isolation, each node $x_k$, $1 \leq k \leq n$, functions as a discrete nonlinear map $f_k$, changing at each iteration $t \in \mathbb{N}$ as $x_k(t) \to x_k(t+1)$. When coupled as a network with adjacency $A$, each node will also receive contributions through the incoming edges from the adjacent nodes. Throughout this paper, we will consider an additive rule of combining these contributions, for a couple of reasons: first, summing weighted incoming inputs is one simple, yet mathematically nontrivial way to introduce the cross talk between nodes; second, summing weighted inputs inside a nonlinear integrating function is reminiscent of certain mechanisms studied in the natural sciences (such as the integrate and fire neural mechanism studied in our previous work in the context of coupled dynamics). The coupled system will then have the following general form:

\begin{eqnarray}
x_k(t) \longrightarrow x_k (t+1) &=& f_k\left(\sum_{k=1}^{n}{g_{jk} A_{jk} x_k} \right) \nonumber 
\label{mothermap}
\end{eqnarray}

\noindent where $g_{jk}$ are the weights along the adjacency edges. One may view this system simply as an iteration of an $n$-dimensional map $f=(f_k)_{k=1}^n$, with $f \colon \mathbb{R}^n \to \mathbb{R}^n$ (in the case of real-valued nodes), or respectively $f \colon \mathbb{C}^n \to \mathbb{C}^n$ (in the case of complex-valued nodes). The new and exciting aspect that we are proposing in our work is to study the dependence of the coupled dynamics on the parameters, in particular on the coupling scheme (adjacency matrix) -- viewed itself as a system parameter. To fix these ideas, we focused first on defining these questions and proposing hypotheses for the case of quadratic node-dynamics. The logistic family is one of the most studied family of maps in the context of both real and complex dynamics of a single variable. It was also the subject of our previous modeling work on random iterations.

In this paper in particular, we will work with quadratic node-maps, with their traditional parametrization $f_c(z) = z^2+c$, with $f_c \colon \mathbb{C} \to \mathbb{C}$ and $c \in \mathbb{C}$ for the complex case and $f_c \colon \mathbb{R} \to \mathbb{R}$ and $c \in \mathbb{R}$ for the real case. The network variable will be called respectively  $(z_1,...z_k) \in \mathbb{C}^n$ in the case of complex nodes, and $(x_1,...x_k) \in \mathbb{R}^n$  in the case of real nodes. We consider both the particular case of identical quadratic maps (equal $c$ values), and the general case of different maps attached to the nodes throughout the network. In both cases, we aim to study the asymptotic behavior of iterated node-wise orbits, as well as of the $n$-dimensional orbits (which we will call multi-orbits). As in the classical theory of Fatou and Julia, we will investigate when orbits escape to infinity or remain bounded, and how much of this information is encoded in the critical multi-orbit of the system. 

For the following definitions, fix the network (i.e., fix the adjacency $A$ and the edge weights $g$). To avoid redundancy, we give definitions for the complex case, but they can be formulated similarly for real maps:

\begin{defn}
For a fixed parameter $(c_1,...,c_n) \in \mathbb{C}^n$, we call the \textbf{filled multi-Julia set} of the network, the locus of $(z_1,...,z_n) \in \mathbb{C}^n$ which produce a bounded multi-orbit in $\mathbb{C}^n$.  We call the \textbf{filled uni-Julia set} the locus of $z \in \mathbb{C}$ so that $(z,...z) \in \mathbb{C}^n$ produces a bounded multi-orbit. The \textbf{multi-Julia set (or the multi-J set)} of the network is defined as the boundary in $\mathbb{C}^n$ of the filled multi-Julia set. Similarly, one defines the \textbf{uni-Julia set (or uni-J set)} of the network as the boundary in $\mathbb{C}$ of its filled counterparts.
\end{defn}

\begin{defn}
We define the \textbf{multi-Mandelbrot set (or the multi-M set)} of the network the parameter locus of $(c_1,...,c_n) \in \mathbb{C}^n$ for which the multi-orbit of the critical point $(0,...,0)$ is bounded in $\mathbb{C}^n$. We call the \textbf{equi-Mandelbrot set (or the equi-M set)} of the network, the locus of $c \in \mathbb{C}$ for which the critical multi-orbit is bounded for \textbf{equi-parameter} $(c_1,c_2,...c_n)=(c,c,...c) \in \mathbb{C}^n$.  We call the \textbf{$k$th node equi-M set} the  locus $c \in \mathbb{C}$ such that the component of the multi-orbit of $(0,...,0)$ corresponding to the $k$th node remains bounded in $\mathbb{C}$.
\end{defn}

\noindent We study, using a combination of analytical and numerical methods, how the structure of the Julia and Mandelbrot sets varies under perturbations of the node-wise dynamics (i.e., under changes of the quadratic multi- parameter $(c_1, c_2, \hdots c_3)$) and under perturbations of the coupling scheme (i.e., of the adjacency matrix $A$ and of the coupling weights $g$). In this paper, we start with investigating these questions in small (3-dimensional) networks, with specific adjacency configurations. In a subsequent paper, we will move to investigate how similar phenomena may be quantified and studied analytically and numerically in high-dimensional networks. In both cases, we are interested in particular in observing differences in the effects on dynamics of three different aspects of the network architecture: (1) increasing/decreasing edge weights, (2) increasing/decreasing edge density, (3) altering edge configuration by adding, deleting or moving edges.

\begin{figure}[h!]
\begin{center}
\includegraphics[width=0.75\textwidth]{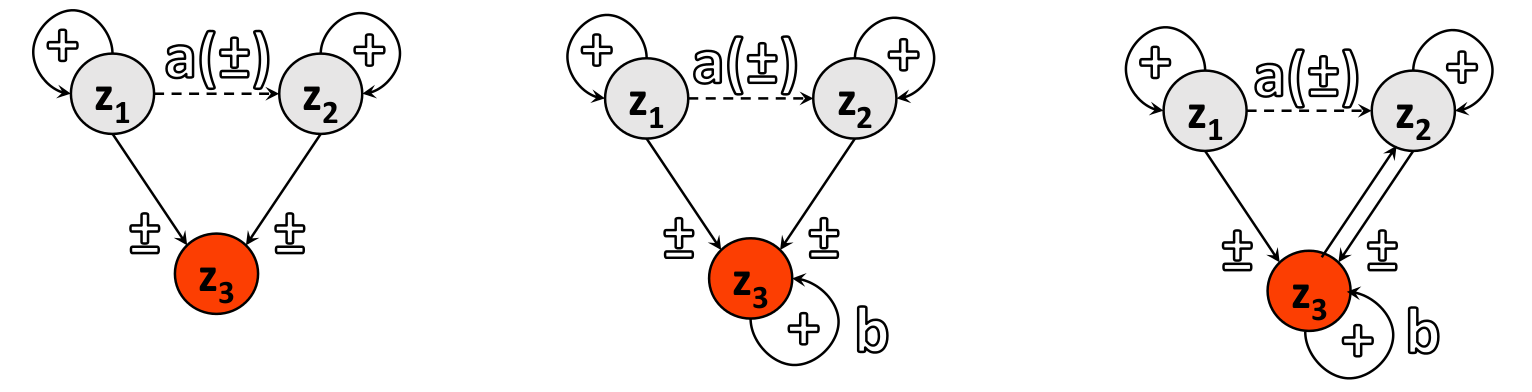}
\end{center}
\caption{\small \emph{{\bf Three dimensional networks} used as simple coupling setups to study the dependence of the Mandelbrot set topology on coupling strength and on network architecture. We will call these three constructions: {\bf A}. the simple dual model; {\bf B.} the self drive model; {\bf C.} the feedback model.}}
\label{3D_net}
\end{figure}

While a desired objective would be to obtain general results for all network sizes (since many natural networks are large), we start by studying simple, low dimensional systems. In this study, we focus on simple networks formed of three nodes, connected by different network geometries and edge weights. To fix our ideas, we will follow and  illustrate three structures in particular (also see Figure~\ref{3D_net}): (1) Two input nodes $z_1$ and $z_2$ are self driven by quadratic maps, and the single output node $z_3$ is driven symmetrically by the two input nodes; $z_1$ additionally communicates with $z_2$ via an edge of variable weight $a$, which can take both positive and negative values. We will call this the \textbf{\emph{simple dual model}}. (2) In addition to the simple dual scheme, the output node $z_3$ is also self-driven, i.e. there is a self-loop on $z_3$ of weight $b$ (which can be positive or negative). We will call this the \textbf{\emph{self-drive model}}. (3) In addition to the self-driven model, there is also feedback from the output node $z_3$ into the node $z_2$, via a new edge of variable weight $f$. We will call this the \textbf{\emph{feedback model}}. Unless specified, edges have positive unit weight. Notice that the same effect as negative feed-forward edges from $z_1$ and $z_2$ into $z_3$ can be obtained by changing the sign of $b$, etc. The three connectivity models we chose to study and compare are described by the equations below:\\

\noindent
\begin{minipage}{0.25\textwidth}
\textbf{\emph{ Simple dual model:}}
\begin{eqnarray}
z_1 &\to& z_1^2 + c_1 \nonumber \\
z_2 &\to& (a z_1 + z_2)^2 + c_2 \nonumber \\
z_3 &\to& (z_1 + z_2)^2 + c_3 \nonumber
\end{eqnarray}
\end{minipage}
\quad \quad \quad
\begin{minipage}{0.3\textwidth}
\textbf{\emph{Self-drive model:}}
\begin{eqnarray}
z_1 &\to& z_1^2 + c_1 \nonumber \\
z_2 &\to& (a z_1 + z_2)^2 + c_2 \nonumber \\
z_3 &\to& (z_1 + z_2 + b z_3)^2 + c_3 \nonumber 
\end{eqnarray}
\end{minipage}
\quad \quad \quad
\begin{minipage}{0.3\textwidth}
\textbf{\emph{Feedback model:}}
\begin{eqnarray}
z_1 &\to& z_1^2 + c_1 \nonumber \\
z_2 &\to& (a z_1 + z_2 + f z_3)^2 + c_2 \nonumber \\
z_3 &\to& (z_1 + z_2 + b z_3)^2 + c_3 \nonumber
\end{eqnarray}
\end{minipage}

\vspace{5mm}
\noindent For a fixed multi-parameter $(c_1,c_2,c_3) \in \mathbb{C}^3$ for example, one can see all three systems as generated by a network map $f=(f_{c_1}$, $f_{c_2}$, $f_{c_3}) \colon \mathbb{C}^3 \to \mathbb{C}^3$, defined as $f(z_1, z_2, z_3) = (f_{c_1}([Az]_1),f_{c_2}([Az]_2),f_{c_3}([Az]_3))$, for any $z=(z_1, z_2, z_3)^t \in \mathbb{C}^3$. 

We try to classify and understand the effects that coupling changes have on the topology of multi-J and multi-M sets for both complex and real networked maps. We don't expect all classical topology results on the Julia and Mandelbrot sets for single maps (e.g., Fatou-Julia theorem, or connectivity of the Mandelbrot set) to carry out for networks of coupled maps. However, since the topology of the full sets in $\mathbb{C}^3$ is somewhat harder to inspect, we study as a first step their equi-slices and node-wise equi-slices, which are objects in $\mathbb{C}$. 

We will track and compare in particular the differences between the three models, but also the geometric and topological changes produced on the equi-slices within each one model for different values of the parameters $a$, $b$ and $f$. None of these results, however, can be directly extrapolated to similar conclusions on the full sets. To offer some insight into the latter, we study the multi-M and multi-J sets in the context of real maps, for which there objects can be visualized in $\mathbb{R}^3$.

\section{Complex coupled maps}
\label{complex_maps}

\subsection{Equi-Mandelbrot sets}
\label{equiM_sets}


A first intuitive question is when the nodes of the network have similar behavior, and whether if one node-wise orbit is bounded, the others will remain bounded. This relationship is trivial to establish in some cases, such as for example in the simple dual model with independent input nodes (i.e., $a=0$). Indeed, in this model, for any fixed $c \in \mathbb{C}$, the origin's orbit in $\mathbb{C}^3$ under $(f_c,f_c,f_c)$ can be described as:
\begin{eqnarray}
z_1&:& 0 \to c \to c^2+c \to (c^2+c)^2 + c \to \hdots \nonumber \\
z_2&:&  0 \to c \to c^2+c \to (c^2+c)^2 + c \to \hdots \nonumber \\
z_3&:&  0 \to c \to (2c)^2 + c \to 4(c^2+c)^2 + c \hdots \nonumber
\end{eqnarray}

\noindent The projection of the orbit in any of the three components only depends on the previous states of $z_1$ and $z_2$, and these three  sequences are simultaneously bounded in $\mathbb{C}$, hence the node-specific equi-Mandelbrot sets are all identical with the traditional Mandelbrot set. Some basic connections between node-wise equi-M sets in each of the three models are stated below. We will prove these incrementally (recall that the dual model is a particular case of self-drive for $b=0$, and the self-drive is a particular case of feedback model with $f=0$).

\begin{figure}[h!]
\begin{center}
\includegraphics[width=\textwidth]{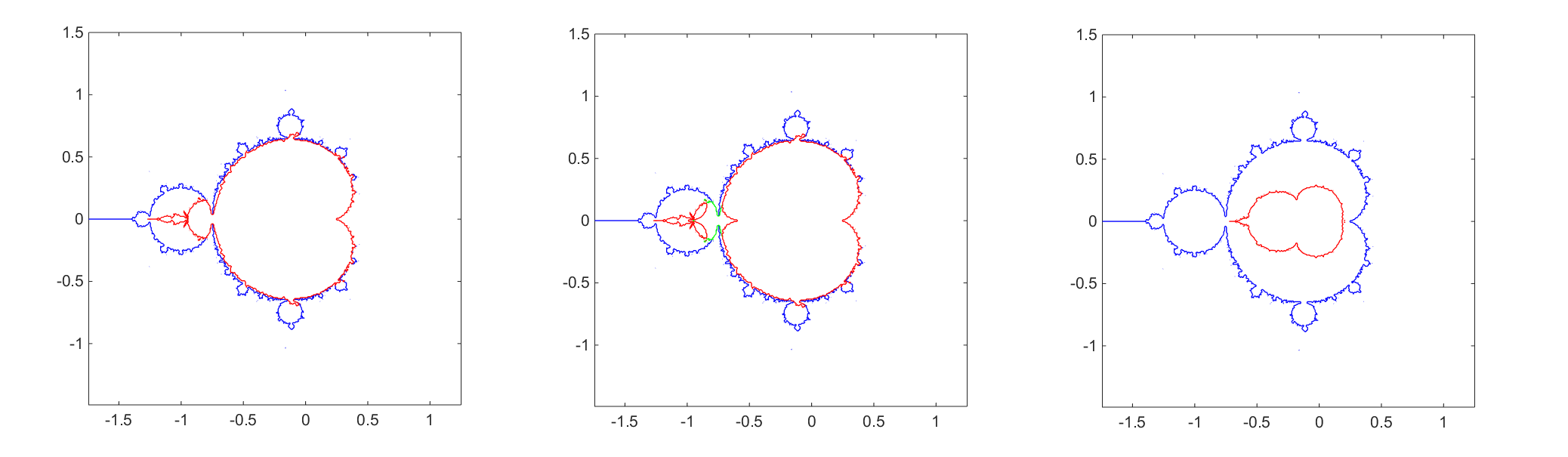}
\end{center}
\caption{\small \emph{{\bf Differences between node-specific equi-Mandelbrot slices}, for different connectivity patterns. {\bf A.} For the simple dual model with $a=-2/3$, the equi-Mandelbrot set for the nodes $z_2$ and $z_3$ are identical (red), but different from the set for the node $z_1$ (blue) {\bf B.} For the self-drive model with negative feedback, $a=-2/3$ and $b=1/3$, the equi-Mandelbrot sets for the three  nodes $z_1$, $z_2$ and $z_3$ (shown respectively in blue, green and red) are all different. {\bf C.} For the feedback model with with negative feedback, $a=-2/3$ and $b=1/3$, $f=-1$, the equi-Mandelbrot set for the nodes $z_2$ and $z_3$ are identical (red), but different from the set for the node $z_1$ (blue). In all panels, the computations were generated based on $L=100$ iterations, and for a test radius of $R=10$.}}
\label{node_differences}
\end{figure}

\begin{prop} In the simple dual model, the node-wise equi-M sets for the nodes $z_2$ and $z_3$ are identical subsets of the traditional Mandelbrot set (which is the equi-M set for node $z_1$). 
\label{prop_simple_dual}
\end{prop}

\proof{
The simple case $a=0$ was already discussed. We will now assume $a \neq 0$. Suppose the critical orbit for node $z_2$ is bounded  by a radius $M$, that is $z_2(n) \leq M$, for all $n$. Hence (omitting the subscript $n$ for simplicity):
$$M \geq \lvert z_2(n+1) \rvert = \lvert (az_1+z_2)^2 + c \rvert \geq \lvert az_1 +z_2 \rvert^2 - \lvert c \rvert \; \Longrightarrow \;  \sqrt{M + \lvert c \rvert} \geq \lvert az_1+z_2 \rvert \geq \lvert az_1 \rvert - \lvert z_2 \rvert$$

\noindent It follows that: $$\lvert az_1 \rvert \leq \sqrt{M + \lvert c \rvert} + \lvert z_2 \rvert \leq \sqrt{M + \lvert c \rvert} + M, \text{ where } a \neq 0$$
Hence if the orbit $z_2$ is bounded, then the orbit of $z_1$ is bounded. This applies in particular for the critical orbit, showing that the equi-M set for $z_2$ if a subset of the equi-M set for $z_1$. 

We will next show that, for the simple dual model, corresponding orbits of $z_2$ and $z_3$ are simultaneously bounded. For instance, suppose that an orbit $z_2(n)$ is bounded by $M >0$. It follows, as shown above, that the corresponding $z_1(n)$ orbit is bounded by some $K>0$. Then:
\begin{eqnarray*}
z_3(n+1) = \lvert (z_1 + z_2)^2 + c \rvert &\leq& \lvert z_1 + z_2 \rvert^2 + \lvert c \rvert = \lvert (az_1+z_2) + (1-a)z_1 \rvert^2 + \lvert c \rvert \\
&\leq& \left( \lvert az_1+z_2 \rvert + \lvert 1-a \rvert \lvert z_1 \rvert \right)^2 + \lvert c \rvert \leq (M+\lvert 1-a \rvert K)^2 + \lvert c \rvert
\end{eqnarray*}

\noindent Hence the orbit $z_3(n)$ is bounded. The converse is similar, showing that the $z_2$ and $z_3$ equi-M sets are always identical subsets of the $z_1$ equi-M set in the simple dual model. In Figure~\ref{node_differences}a, we show that these are generally strict subsets, and that a non-symmetric communication $a \neq 0$  can introduce significant differences between the traditional  Mandelbrot set of the independent node $z_1$ and the equi-M subsets for $z_2$ and $z_3$. For example, it is not hard to show that, for $a=-2/3$ (illustrated in Figure~\ref{node_differences}a), the point $c=-2$ belongs to the Mandelbrot set of $z_1$ (the critical orbit has period three), but not to the equi-M set of $z_2$ and $z_3$. \\
\qed}

\vspace{3mm}
\noindent An additional self-drive $b \neq 0$ applied to the output node changes the balance of inputs to $z_3$, in the following sense:

\begin{prop} In the self-drive model, the  node-wise equi-M sets of $z_2$ and $z_3$ remain subsets of the standard Mandelbrot set, but the equi-M set of $z_3$ is strictly contained in the equi-M set of $z_2$ (Figure~\ref{node_differences}b). 
\label{prop_self_drive}
\end{prop}

\proof{To prove the first part of this statement, take a point $c$ in the equi-M set of $z_3$, meaning that the orbit of $z_3$ is bounded: there exists $M>0$ such that $\lvert z_3(n) \rvert \leq M$, for all $n \geq 0$. We can express:
$$\lvert z_3(n+1) \rvert = \lvert (z_1 + z_2 +bz_3)^2 + c \rvert \geq \lvert z_1 + z_2 +bz_3 \rvert^2 - \lvert c \rvert$$

\noindent It follows that:
$$\lvert z_1 + z_2 \rvert - \lvert bz_3 \rvert \leq \lvert z_1+z_2 +bz_3 \rvert \leq \sqrt{\lvert z_3(n+1) \rvert + \lvert c \rvert} \leq \sqrt{M + \lvert c \rvert}$$

\noindent Hence 
$$\lvert z_1 + z_2 \rvert \leq \sqrt{M + \lvert c \rvert}  +  \lvert bz_3 \rvert \leq \sqrt{M + \lvert c \rvert} + \lvert b \rvert M$$

\noindent that is, the sequence $\xi(n) = z_1(n) + z_2(n)$ is also bounded in radius by $K_1 =  \sqrt{M + \lvert c \rvert} + \lvert b \rvert M$. Let us recall that 
\begin{eqnarray*}
\lvert \xi(n+1) \rvert &=& \lvert z_1^2 + c + (az_1+z_2)^2 + c \rvert = \lvert z_1^2 +  [(a-1)z_1 + \xi]^2 + 2c \rvert\\\\
&=& z_1^2 + (a-1)^2z_1^2 + 2(a-1)z_1 \xi + \xi^2 +2c \\\\
&=& \left \lvert  \left[ \sqrt{(a-1)^2+1}z_1 + \frac{(a-1)\xi}{\sqrt{(a-1)^2+1}} \right]^2 + \xi^2 \left(1- \frac{(a-1)^2}{(a-1)^2+1} \right) + 2c \right \rvert\\\\
&\geq& \left \lvert  \sqrt{(a-1)^2+1}z_1 + \frac{(a-1)\xi}{\sqrt{(a-1)^2+1}} \right \rvert^2 -  \left \lvert  \frac{\xi^2}{(a-1)^2+1}  + 2c \right \rvert
\end{eqnarray*}

\noindent It follows that
$$\lvert k_1z_1+k_2\xi \rvert^2 \leq \lvert \xi(n+1) \rvert + \left \lvert  \frac{\xi^2}{(a-1)^2+1}  + 2c  \right \rvert \; \Longrightarrow \; \lvert k_1z_1+k_2\xi \rvert \leq K_2$$ 

\noindent where $k_1 = \sqrt{(a-1)^2+1}$, $\ds k_2 = \frac{(a-1)}{\sqrt{(a-1)^2+1}}$ and $\ds K_2 = \sqrt{M + \frac{M}{(a-1)^2+1} + \lvert 2c \rvert}$. It follows that $\lvert k_1 z_1 \rvert \leq K_2 + \lvert k_2 \xi \rvert \leq K + \lvert k_2 \rvert M$, hence the orbit of the node $z_1$ is also bounded. Now recall that: $\xi = z_1 + z_2$ is bounded. Since we can write $\lvert z_2 \rvert = \lvert \xi - z_1 \rvert \leq \lvert \xi \rvert + \lvert z_1 \rvert$, it follows that $z_2$ is also bounded. This proves that the equi-M set of $z_3$ is a subset of the equi-M set of $z_2$, which is in turn a subset of the traditional Mandelbrot set (i.e., the equi-M set of $z_1$).

To prove that these inclusions are strict, one can easily find points which are in the equi-M set of $z_2$, but not in the equi-M set of $z_3$. For example, for the parameters in Figure Figure~\ref{node_differences}b, $c=-3/4$ is in the M-set of $z_1$ (the critical orbit is eventually fixed), and it is in the equi-M set of $z_2$, but it in not in the M-set of $z_3$. Indeed, for these particular parameters, $z_2(n+1) = z_2^2(n) + z_2(n) -1/2$, with $z_2(2) = -3/4$. One can easily show that, if $z_2 \in [-1,0]$, then $z_2^2 + z_2 -1/2 \in [-1,0]$, hence it follows by induction that the critical orbit of $z_2$ is contained in $[-1,0]$ (i.e., bounded).\\
\qed}

\vspace{3mm}
\noindent Finally, introducing any arbitrary feedback $f \neq 0$ re-couples the behavior of nodes $z_2$ and $z_3$, producing a common equi-Mandelbrod set, largely shrunk from the simple dual version:

\begin{prop} In the feedback model with $b \neq 0$ and $f \neq 0$, the node-wise equi-M sets for the nodes $z_2$ and $z_3$ are again identical subsets of the traditional Mandelbrot set (Figure~\ref{node_differences}c).
\label{prop_feedback}
\end{prop}

\proof{The proof is a slightly more general version of that for Proposition~\ref{prop_self_drive}. Suppose first that the orbit of $z_3$ is bounded  in radius by $M$. As before, it follows that:
$$\ds M \geq \lvert z_3(n+1) \rvert = \lvert (z_1+z_2 + bz_3)^2 + c \rvert \geq \lvert z_1+z_2 +bz_3 \rvert^2 - \lvert c \rvert$$

\noindent Hence, as before $$\lvert \xi \rvert = \lvert z_1 + z_2 \rvert \leq \lvert b \rvert M + \sqrt{M+\lvert c \rvert} = K_1$$

\noindent Call $\psi = \xi + fz_3$ so that:
$$\lvert \psi \rvert \leq \lvert \xi \rvert + \lvert fz_3 \rvert \leq K_1 + \lvert f \rvert M = K_2$$

\noindent We calculate:
\begin{eqnarray*}
\lvert \xi(n+1) \rvert &=& \lvert z_1^2 + (az_1+z_2+fz_3)^2 + 2c  \rvert = \lvert z_1^2 + [(a-1)z_1+\xi+fz_3)^2 + 2c  \rvert\\ \\
&=&  \lvert z_1^2 + [(a-1)z_1+ \psi]^2 + 2c  \rvert = \lvert z_1^2 + (a-1)^2z_1^2 + 2(a-1)z_1\psi +\psi^2 + 2c \rvert \\\\
&=& \left \lvert \left( [1+ (a-1)^2]z_1^2 + 2(a-1)z_1 \psi + \frac{a-1}{(a-1)^2+1}\psi^2  \right) + \left( 1- \frac{a-1}{(a-1)^2+1}  \right) \psi^2  + 2c \right \rvert\\\\
&=& \left \lvert (k_1z_1 + k_2 \psi)^2 + k_3\psi^2 + 2c \right \rvert
\end{eqnarray*}

\noindent where $k_1 = \sqrt{(a-1)^2+1}$, $k_2 = \ds \frac{a-1}{(a-1)^2+1}$ and $\ds k_3 = \frac{1}{(a-1)^2+1}$. Hence
$$\lvert k_1 z_1 + k_2 \psi \rvert \leq \lvert \xi(n+1) \rvert + k_3 \lvert \psi \rvert^2 + 2\lvert c \rvert \leq K_1 + k_3 K_2^2 + 2\lvert c \rvert = K_3$$

\noindent Since $\psi$ is bounded, it follows that $z_1$ is bounded. Since $\xi$ is bounded, it follows that $z_2$ is bounded. This proofs that the equi-M set of $z_3$ is a subset of the equi-M sets of $z_2$ and $z_1$.\\

\noindent Conversely, suppose that the orbit of $z_2$ is bounded in radius by $M>0$. It follows that the sequence corresponding to $az_1 + fz_3$ is also bounded. Recall that the case $f=0$ was covered by the previous proposition, hence we can assume now that $f \neq 0$. Hence $\xi = hz_1 + bz_3$ is bounded by a constant $K_1$, where $h=ba/f$, for $f \neq 0$. As before, call $\psi = hz_1 + bz_3 + z_2$, and notice that $\psi$ is bounded by $K_1+ M$.
\begin{eqnarray}
\lvert \xi(n+1) \rvert &=& \lvert h(z_1^2 + c) + b[(z_1+z_2+bz_3)^2 + c] \rvert = \lvert hz_1^2 + b[\psi + (1-h)z_1]^2 + c(h+b) \rvert \nonumber \\ \nonumber \label{eqn1}\\
&=& \lvert  hz_1^2 + b(1-h)^2z_1^2 + 2b(1-h)z_1\psi + b\psi^2  +c(h+b) \rvert 
\end{eqnarray}

\noindent If $\ds k_1^2 = h + b(1-h)^2 \neq 0$, then we again have:
\begin{eqnarray*}
\lvert \xi(n+1) \rvert &=& \lvert (k_1z_1 + k_2 \psi)^2 + k_3 \psi^2 + c(h+b) \rvert
\end{eqnarray*}

\noindent where $\ds k_2 = \frac{1-h}{k_1}$ and $k_3 = f - k_2^2$. Since $\psi$ is bounded, it follows that $z_1$ is bounded. Since $\xi$ is bounded, it further follows that $z_3$ is bounded.

We look separately at the case  $\ds k_1^2 = h + b(1-h)^2 = 0$, for which Equation~\eqref{eqn1} becomes:
\begin{eqnarray*}
\lvert \xi(n+1) \rvert &=& \lvert 2b(1-h)z_1\psi + b\psi^2  +c(h+b) \rvert 
\end{eqnarray*}

\noindent Since $b \neq 0$, it follows that $b(1-h) \neq 0$. Since $\xi$ and $\psi$ are bounded, it immediately follows that $z_1$ is also bounded.\\

\noindent This concludes the proof that the equi-M sets of $z_2$ and $z_3$ are identical, and both subsets of the equi-M set of $z_1$.
\qed}

\vspace{1cm}
\noindent For the rest of the section, the term of ``equi-M set'' will be referring to the equi-Mandelbrot set of the network, which is the intersection of the three node-specific sets. We illustrate the equi-M set for the three models and for different levels of cross-talk $a$, $b$ and $f$ between nodes. 

\begin{figure}[h!]
\begin{center}
\includegraphics[width=\textwidth]{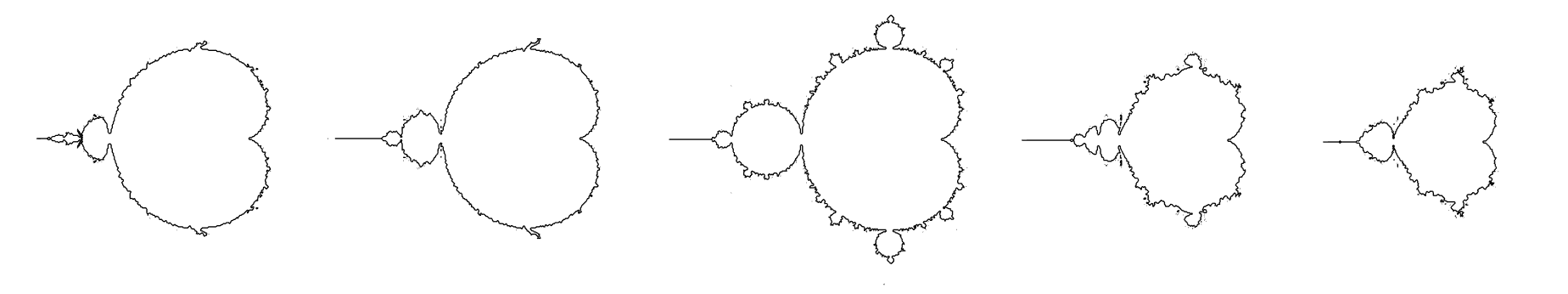}
\end{center}
\caption{\small \emph{{\bf Differences between equi-Mandelbrot slices} in the case of the simple dual model, as the cross-talk parameter $a$ increases: {\bf A.} $a=-2/3$; {\bf B.} $a=-1/3$; {\bf C.} $a=0$ (traditional Mandelbrot set); {\bf D.} $a=1/3$; {\bf E.} $a=2/3$.}}
\label{dual_input}
\end{figure}

Starting with the simple dual input version of the model, we show in Figure~\ref{dual_input} the effects of changing the level $a$ of talk between the input nodes, on the shape of the equi-M set. It is not surprising that, in both positive and negative $a$ ranges, increasing $\lvert a \rvert$ gradually shrinks the equi-Mandelbrot set. This can be motivated intuitively by the fact that an additional contribution to the node $z_2$ may cause the critical orbit to increase faster in the $z_2$, and subsequently the $z_3$ components, hence points in the traditional set will no longer be included in the mutants for $a \neq 0$. 

As $a$ increases in the positive range, we noticed that the network M sets form nested subsets (which is not true for the negative range), that they remain connected for all values of $a$, and that the Hausdorff dimension of the boundary increases with $a$ (in Figure~\ref{dual_input}, notice an increased wrinkling of the boundary as $a$ takes larger positive values, and an increase smoothing as $a$ takes negative values with increasing absolute value). Perturbations of $a$ in the positive range seem to have a much more substantial contribution to the size of the equi-M set, while perturbations of $a$ in the negative range have a lesser influence on the size, and affect mostly the region close tot the boundary of the equi-M set, and the boundary topological details. We will track the same changes in $a$ in the other network models, and investigate if this trend is consistent.

\begin{figure}[h!]
\begin{center}
\includegraphics[width=0.9\textwidth]{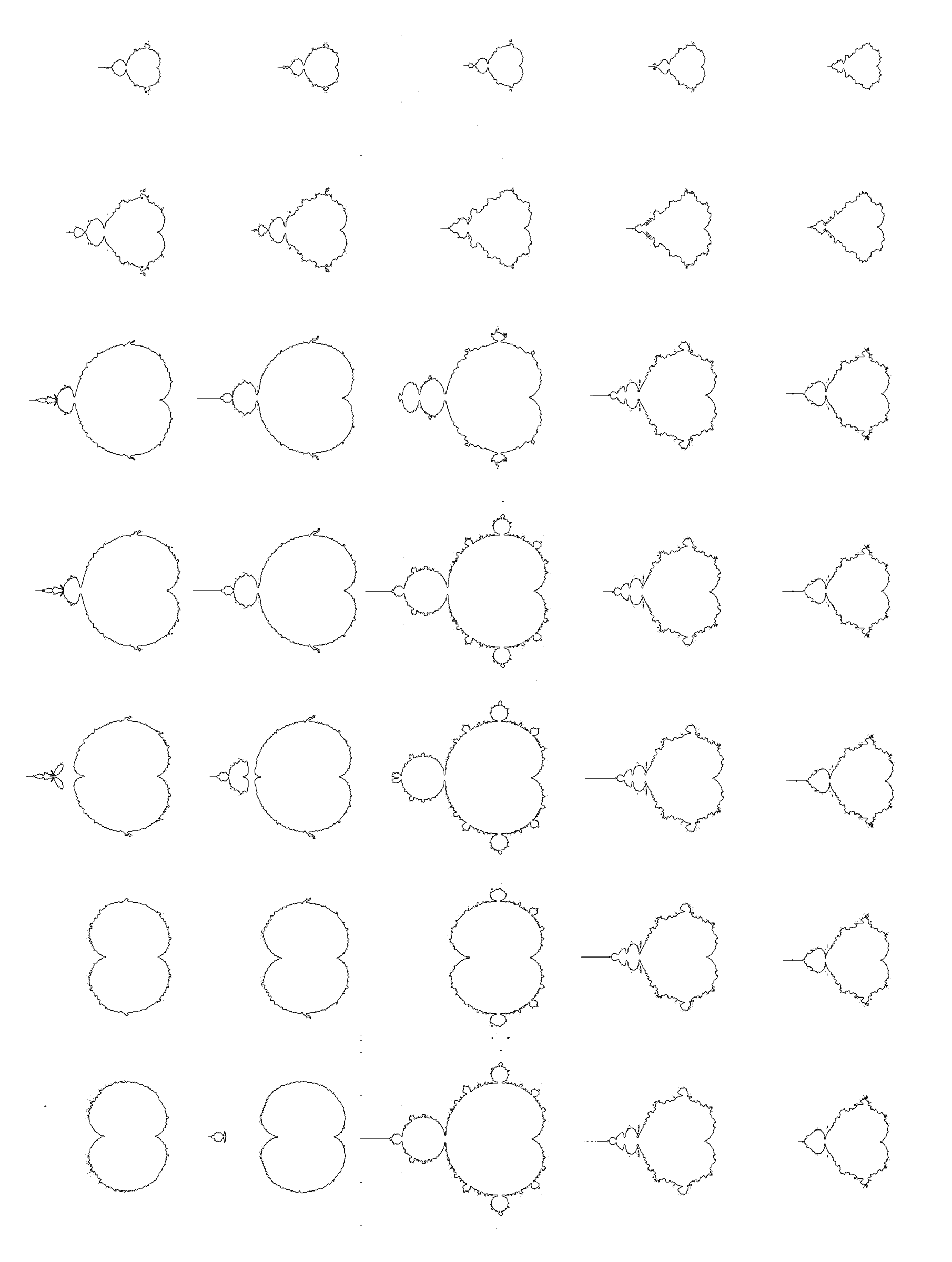}
\end{center}
\caption{\small \emph{{\bf  Equi-M sets for the model variation 2 (with self-drive), for different values of the parameters $a$ and $b$.} The rows show, from bottom to top, increasing values of the self-drive: $b=-1$, $b=-2/3$; $b=-1/3$; $b=0$ (this row representing the simple dual model, as shown in Figure~\ref{dual_input}); $b=1/3$; $b=2/3$; $b=1$. The columns show, from left to right, increasing values of cross-talk between the two input nodes: $a=-2/3$, $a=-1/3$, $a=0$, $a=1/3$ and $a=2/3$. All the equi-M sets were generated from $L=100$ iterations, and plotted at the same scale, in the complex square $[-1.75,1.25] \times [-1.5,1.5]$.}}
\label{equiMand_selfdrive}
\end{figure}

\noindent Figure~\ref{equiMand_selfdrive} illustrates the evolution of the equi-M set in the case of the model with self-drive, for a grid of positive and negative values of the input connectivity $a$ and of the self-drive $b$. Below are some simple visual observations based on our numerical computations, to be addressed analytically in future work.

Decreasing $b$ in the negative range produces no alteration of the M sets when $a>0$. However, it induces dramatic changes in shape and connectivity when $a<0$. If for $b>0$ relatively large, increasing $a$ only slightly alters the shape of the set, for small $b>0$ the size of the set is also altered with increasing $a$ (generating smaller and smaller subsets), and the complexity of its boundary also seems to increase.  The effects of varying $a<0$ for a fixed value of $b$ become more dramatic with decreasing $b$ in the negative range. These effects include changes in shape and topology, the region $a<0$ and $b<0$ allowing the M-set to break into multiple connected components.

\subsection{Uni-Julia sets}
\label{uniJ_sets}

In this section, we will track the changes in the uni-Julia set when the parameters of the system change. One of our goals is to test, first in the case of equi-parameters $c \in \mathbb{C}$, then for general parameters in $\mathbb{C}^3$, if a  Fatou-Julia type theorem applies in the case of our three networks.

First, we try to establish a hypothesis for connectedness of uni-J sets, by addressing numerically and visually questions such as: ``Is it true that if $c$ is in the equi-M set of a network, then the uni-Julia set is connected?'' ``Is it true that, if $c$ is not in the equi-M set of the network, then the uni-Julia set is totally disconnected?'' Clearly, this is not simply a $\mathbb{C}^3$ version of the traditional Fatou-Julia theorem, but rather a slightly different result involving the projection of the Julia set onto a uni-slice. Notice that a connected uni-J set in $\mathbb{C}$ may be obtained from a disconnected $\mathbb{C}^3$ network Julia set, and conversely, that a disconnected uni-J projection may come from a connected Julia in $\mathbb{C}^3$. We will further discuss $\mathbb{C}^3$ versions of these objects in the context of iterations of real variables, where one can visualize the full Mandelbrot and Julia sets for the network as subsets of $\mathbb{R}^3$. Here, we will first investigate uni-J sets for equi-parameters $(c_1,c_2,c_3) = (c,c,c)$, with a particular focus on tracking the topological changes of the uni-J set as the system approaches the boundary of the equi-M set and leaves the equi-M set.

\begin{figure}[h!]
\begin{center}
\includegraphics[width=0.65\textwidth]{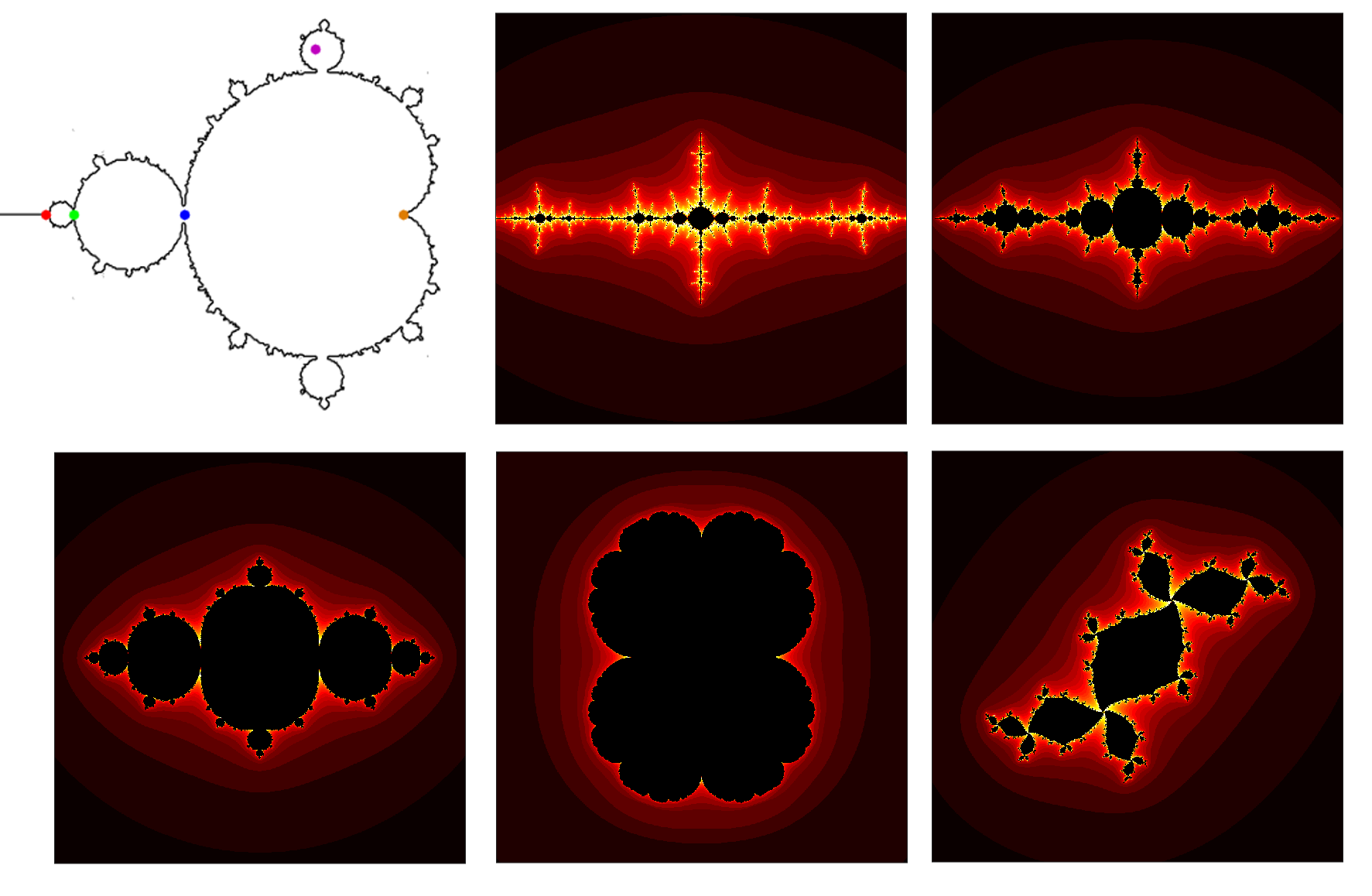}
\end{center}
\caption{\small \emph{{\bf Uni-Julia sets} for a self-drive network with $a=0$ and $b=-1$, for different values of the equi-parameter $c$ (marked with colored dots on the equi-M template in upper left): $c=-1.38$ (red); $c=-1.25$ (green); $c=-0.75$ (blue); $c=0.25$ (orange); $c=-0.15+0.75i$ (purple). All sets were based on $100$ iterations. Both equi-M and uni-J sets coincide in this case with the traditional Mandebrot and Julia sets for single map iterations.}}
\label{uniJulia1}
\end{figure}

\begin{figure}[h!]
\begin{center}
\includegraphics[width=0.9\textwidth]{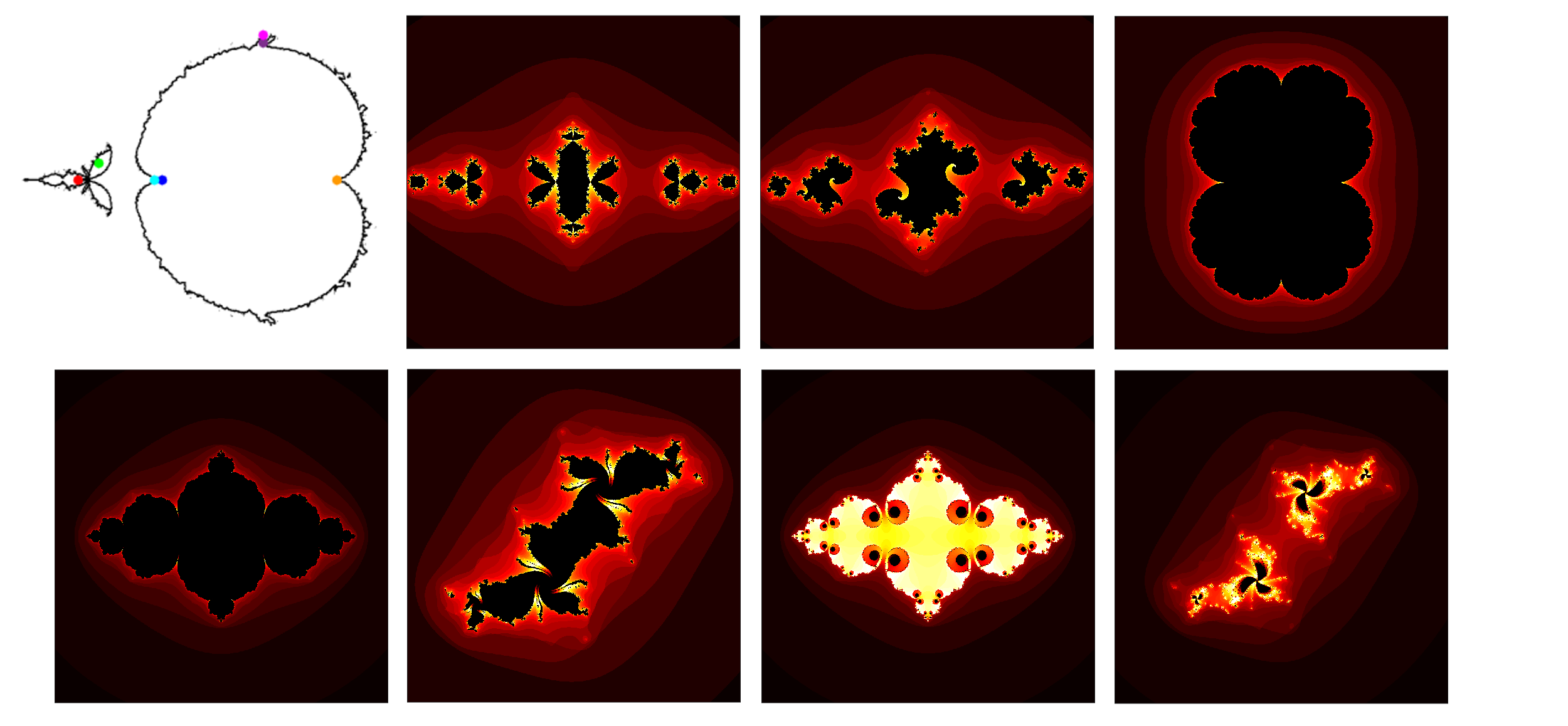}
\end{center}
\caption{\small \emph{{\bf Uni-Julia sets} for a self-drive network with $a=-2/3$ and $b=-1/3$, for different values of the equi-parameter $c$ (marked with colored dots on the equi-M template in upper left): $c=-1$ (red); $c=-0.9+0.08i$ (green); $c=0.25$ (orange); $c=-0.595$ (blue); $c=-0.11+0.66i$ (dark purple); $c=-0.63$ (cyan); $c=-0.11+0.7i$ (magenta). For the first four panels, $c$ is in the equi-M set; for the last two, $c$ is outside of the equi-M set. All sets were based on 100 iterations.}}
\label{uniJulia2}
\end{figure}

\begin{figure}[h!]
\begin{center}
\includegraphics[width=0.6\textwidth]{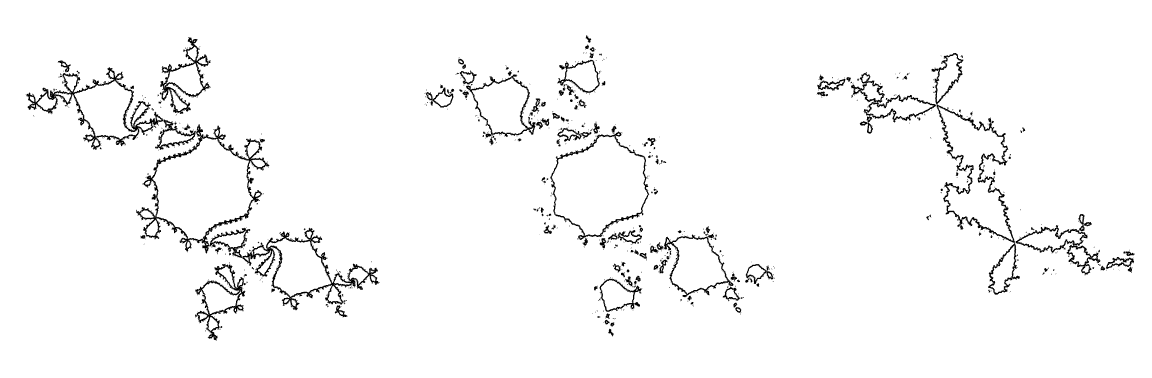}
\end{center}
\caption{\small \emph{{\bf Evolution of the uni-Julia set} for fixed equi-parameter $c=-0.117-0.76i$, as the network profile is changed, from {\bf A.} simple dual with $a=-0.05$, to {\bf B.} self-drive with additional $b=-1$, to {\bf C.} feedback with additional $f=-0.75$.}}
\label{uniJulia_model_comparison}
\end{figure}

First, we fix the network type and the connectivity profile (i.e., the parameters $a$, $b$ and $f$), and we observe how the uni-J sets evolves as the equi-parameter $c$ changes. In Figures~\ref{uniJulia1} and \ref{uniJulia2} we illustrate this for two examples of self-driven models: one with $b=-1$ and $a=0$, the other with $b=-1/3$ and $a=-2/3$. As the parameter point $(c,c,c)$ approaches the boundary of the equi-M set, the topology of the uni-J set if affected, with its connectivity braking down ``around'' the boundary. 

Second, we look at the dependence of uni-Julia sets on the coupling profile (network type). As an example, we fixed the equi-parameter $c=-0.62-0.432i$, and we first considered a simple dual network with negative feed-forward and small cross-talk $a=-0.01$. We then added self-drive $b=-1$ to the output node, then additionally introduced a small negative feedback $f=-0.1$. The three resulting uni-Julia sets are shown in Figure~\ref{uniJulia_model_comparison}. Notice that a very small degree of feedback $f$ produces a more substantial difference than a significant change in the self-drive $b$.

Third, one can study the dependence of uni-Julia sets on the strength of specific connections within the network. As a simple illustration of how complex this dependence may be, we show in Figures~\ref{c2_Julia} and~\ref{c3_Julia} the effects on the uni-J sets of slight increases in the cross-talk parameter $a$, for two different values of the equi-parameter $c$.

An immediate observation is that uni-J sets no not exhibit the dichotomy from traditional single-map iterations no longer stands: uni-J sets can be connected, totally disconnected, but also disconnected into a (finite or infinite) number of connected components, without being totally disconnected. Based on our illustrations, we can further conjecture, in the context of our three models, a description of connectedness for uni-J sets, as follows:

\begin{figure}[h!]
\begin{center}
\includegraphics[width=0.9\textwidth]{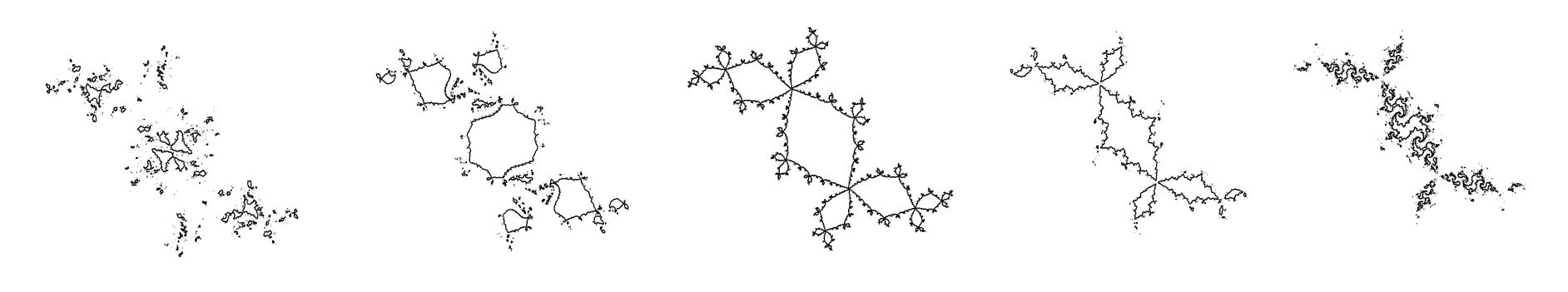}
\end{center}
\caption{\small \emph{{\bf Evolution of the uni-Julia set} for a self-drive network with $b=-1$ and equi-parameter $c=-0.117-0.76$, as the input cross-talk $a$ is increased. The panels show, left to right: $a=-0.07$, $a=-0.05$, $a=0$, $a=0.05$ and $a=0.07$.}}
\label{c2_Julia}
\end{figure}

\begin{figure}[h!]
\begin{center}
\includegraphics[width=0.9\textwidth]{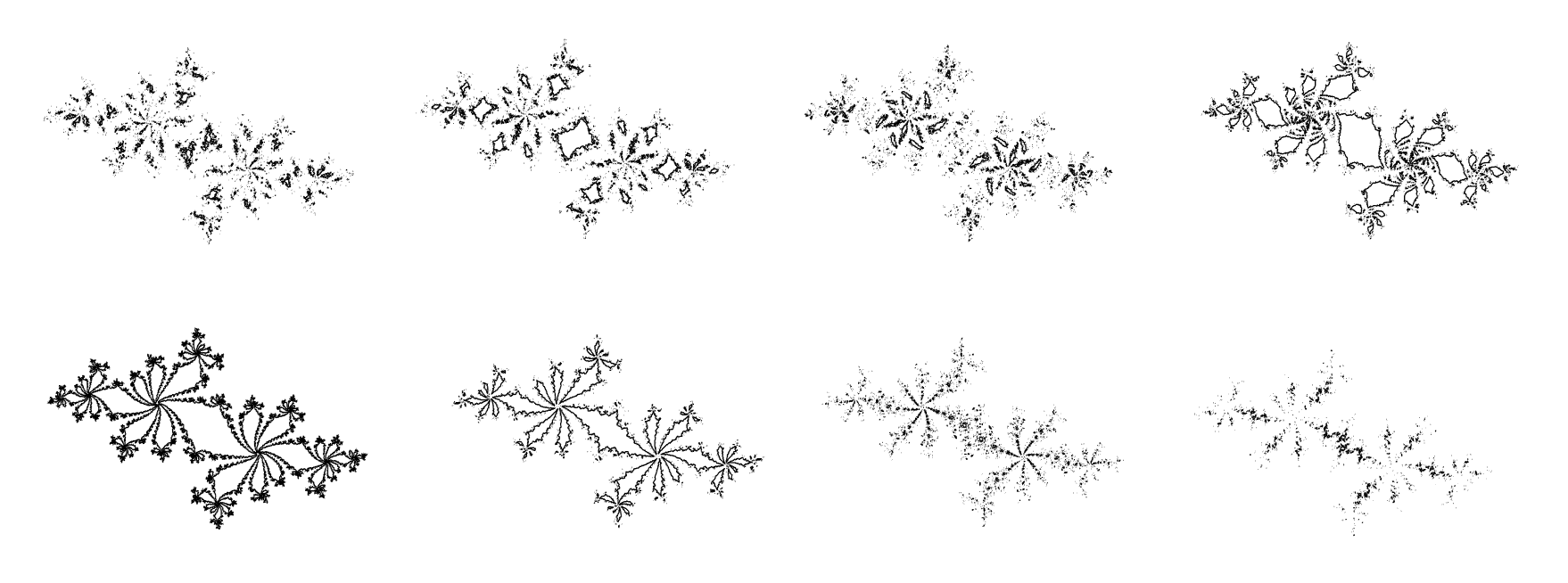}
\end{center}
\caption{\small \emph{{\bf Evolution of the uni-Julia set} for a self-drive network with $b=-1$ and equi-parameter $c=-0.62-0.432i$, as the input cross-talk $a$ is increased. The panels show, left to right: $a=-0.022$, $a=-0.02$,  $a=-0.015$, $a=-0.01$, $a=0$, $a=0.01$ and $a=0.015$, $a=0.02$.}}
\label{c3_Julia}
\end{figure}

\begin{conj} For any of the three models described, and for any equi-parameter $c \in \mathbb{C}$, the uni-J set is connected only if $c$ is in the equi-M set of the network, and it is totally disconnected only if $c$ is not in the equi-M set of the network.
\end{conj}

\noindent {\bf Remark.} The conjecture implies a looser dichotomy regarding connectivity of uni-J sets than that delivered by the traditional Fatou-Julia result for single maps: If $c$ is in the equi-M set of the network, then the uni-J set is either connected or disconnected, without being totally disconnected. If $c$ is not in the equi-M set of the network, then the uni-J set is disconnected (allowing in particular the case of totally disconnected).\\

\begin{figure}[h!]
\begin{center}
\includegraphics[width=0.7\textwidth]{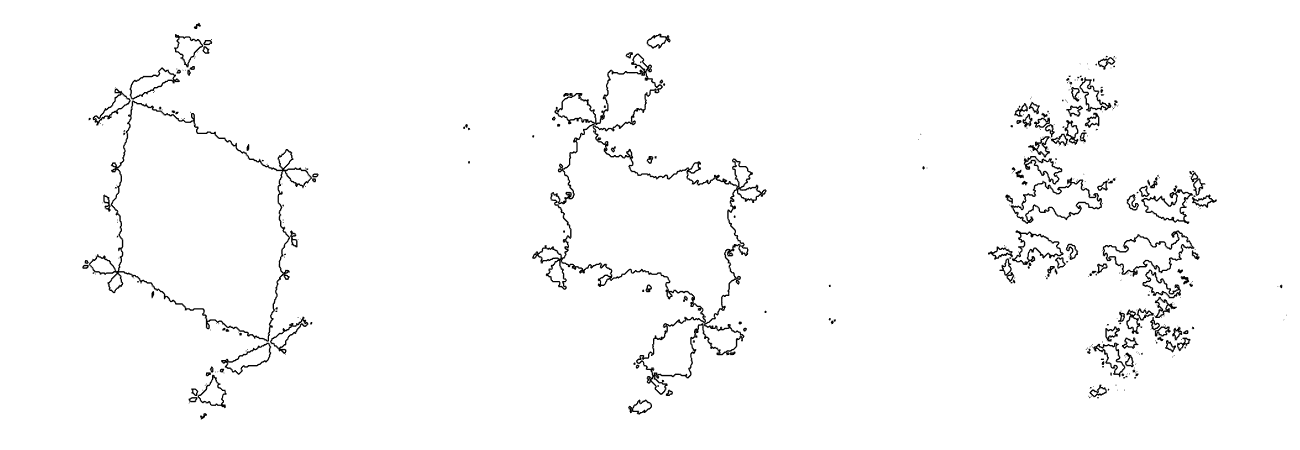}
\end{center}
\caption{\small \emph{{\bf Uni-Julia sets} for a general parameter $(c_1,c_2,c_3) \in \mathbb{C}^3$, with $c_1=-0.75$, $c_2=-0.117-0.76i$ and $c_3=-0.62-0.432$. the panels represent uni-J sets for a self-drive network with $b=-1$, as the cross-talk $a$ changes from {\bf A.} $a=0$, to {\bf B.} $a=0.1$, to {\bf C.} $a=0.15$.}}
\label{general_uniJulia}
\end{figure}

\noindent Finally, we want to remind the reader that uni-Julia sets can be defined for general parameters $(c_1,c_2,c_3) \in \mathbb{C}^3$, as shown in Figure~\ref{general_uniJulia}.

\section{Real case}
\label{real_maps}

\begin{figure}[h!]
\begin{center}
\includegraphics[width=0.9\textwidth]{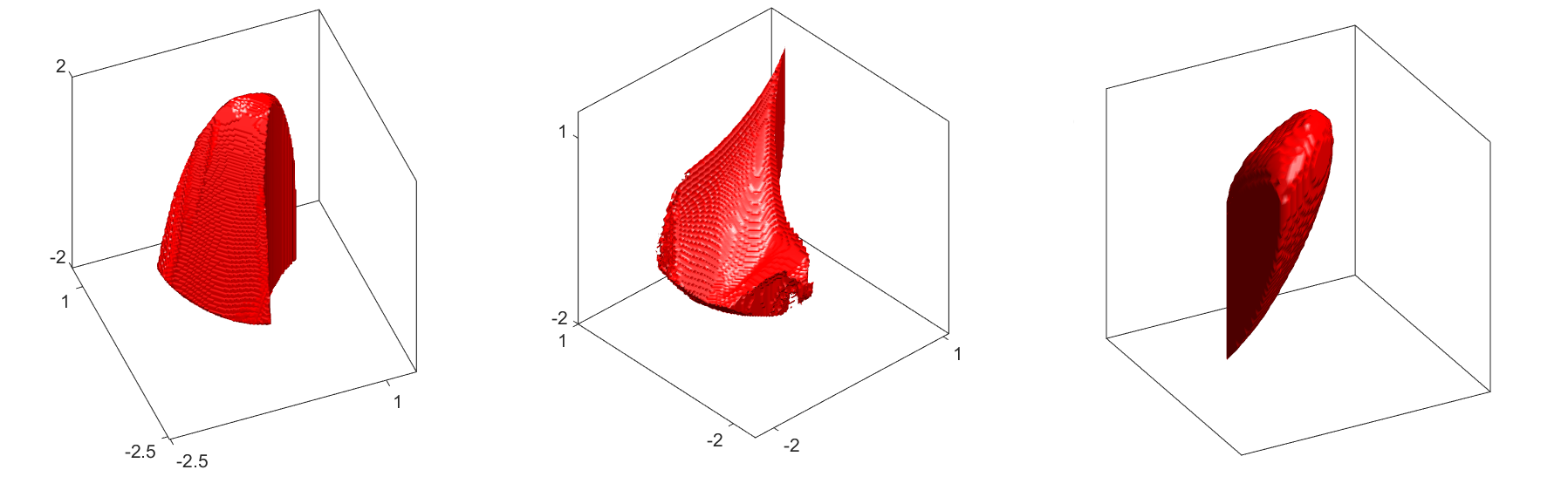}
\end{center}
\caption{\small \emph{{\bf Real network Mandelbrot sets.} {\bf Left.} Simple dual network with $a=-1$. {\bf Center.} Self-drive network with $a=-1$ and $b=1$. {\bf Right.} Self-drive network with $a=1/2$ and $b=1$. Plots were generated with $50$ iterations and in resolution $200^3$.}}
\label{realMand}
\end{figure}

\begin{figure}[h!]
\begin{center}
\includegraphics[width=0.9\textwidth]{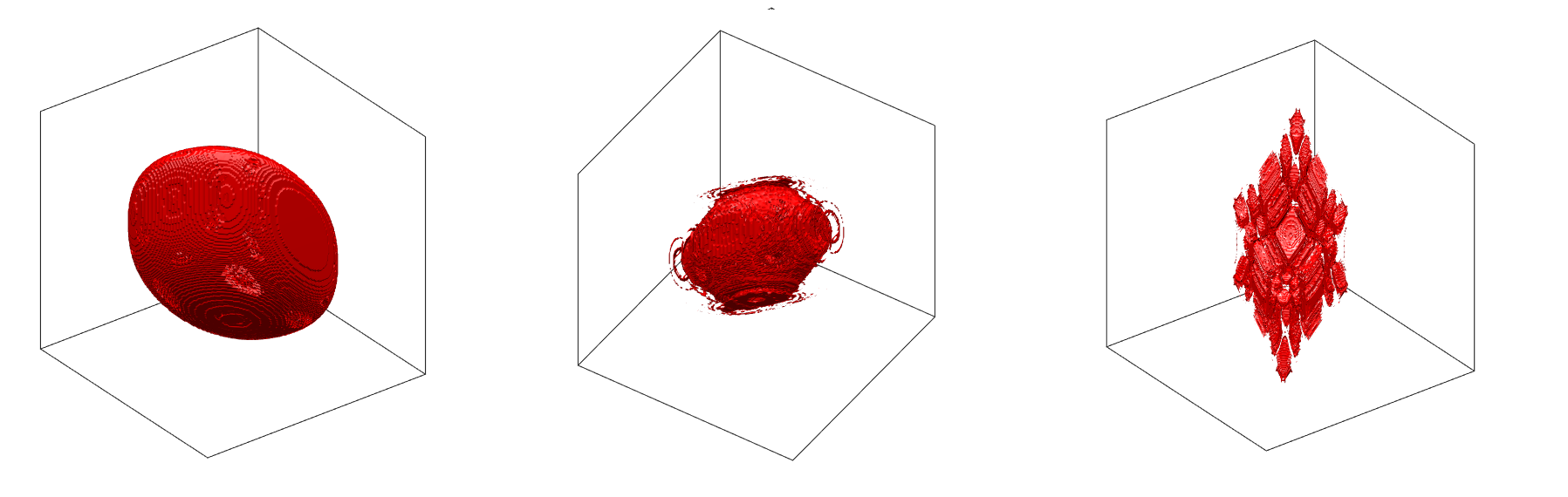}
\end{center}
\caption{\small \emph{{\bf Real network Julia sets.} All panels represent self-drive networks with $a=1/2$ and $b=1$, with equi-parameters respectively: {\bf A.} $c=-0.589$; {\bf B.} $c=-0.4-0.08i$; {\bf C.} $c=-0.62-0.432i$. Plots were generated with $50$ iterations and in resolution $200^3$.}}
\label{realJulia}
\end{figure}

The same definitions apply for iterations of real quadratic maps, with the real case presenting the advantage of easy visualization of full Julia and Mandelbrot sets, rather than having to consider equi-slices, as we did in the complex case. In Figures~\ref{realMand} and~\ref{realJulia}, we illustrate a few multi-M  and multi-J sets respectively, for some of the same networks considered in our complex case. 

Moving to illustrate the \emph{relationship} between the multi-M and the multi-J set in this case, consider for example the self-drive real network with $a=1/2$ and $b=1$, for different parameters $(c_1,c_2,c_3)$. While more computationally intensive, higher-resolution figures would be necessary to establish the geometry and fractality of these sets, one may already notice basic properties. For example, Figure~\ref{realJulia} shows that, if one were to consider complex equi-parameters, the multi-Julia set may not only be connected (Figure~\ref{realJulia}a), or totally disconnected (not shown), but may also be broken into a number of connected components (Figures~\ref{realJulia}b anc c).

\begin{figure}[h!]
\begin{center}
\includegraphics[width=0.6\textwidth]{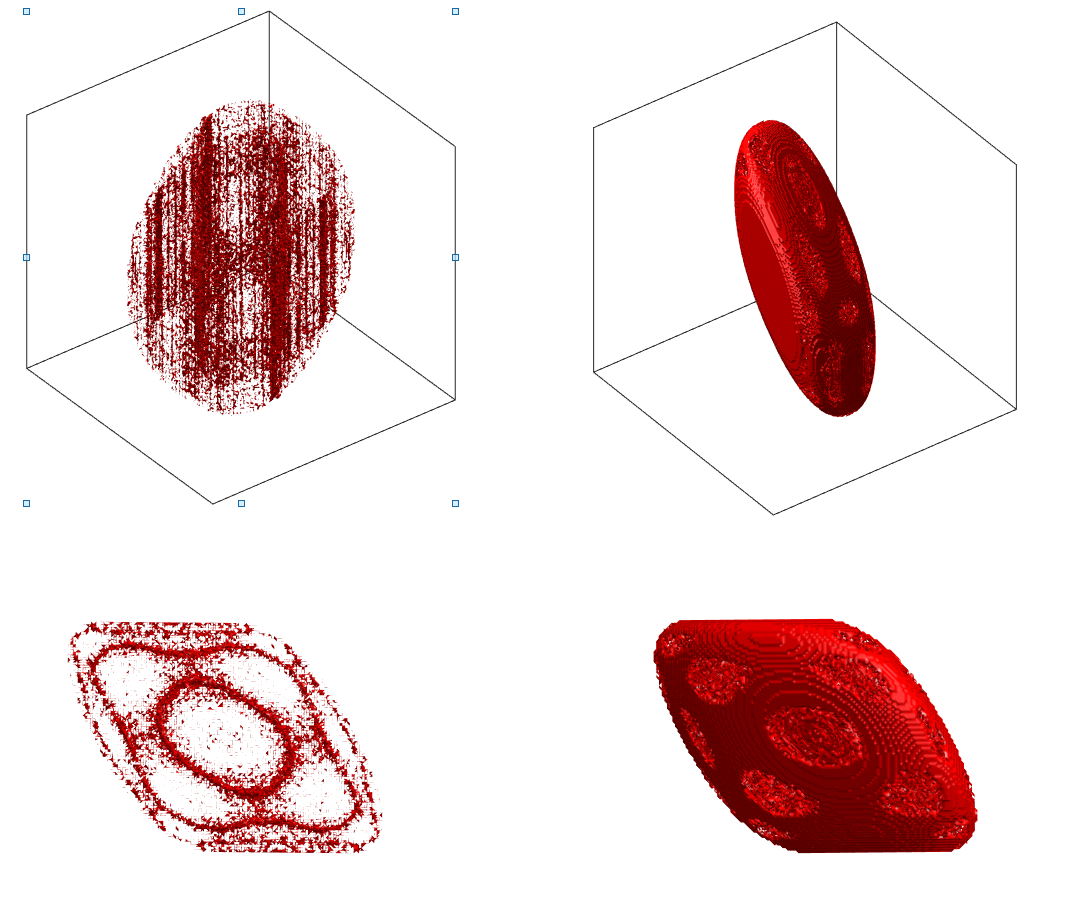}
\end{center}
\caption{\small \emph{{\bf Real network Julia sets} for the self-drive model with $a=1/2$ and $b=1$. The two multi-parameters $(c_1,c_2,c_3)=(-0.5, -0.7, -0.7)$ (left panels) and $(c_1,c_2,c_3)=(-0.5, -0.7, -0.6)$ (right panels) are not in the Mandelbrot set for the network. The top row shows the 3-dimensional Julia sets, the bottom panels show top views of the same sets. Plots were generated with $50$ iterations and in resolution $200^3$.}}
\label{real_comp}
\end{figure}

This remained true if we returned to our restriction of having real parameters, once we allow arbitrary (that is, not necessarily equi) parameters. The panels of Figure~\ref{real_comp} show the multi-J sets for two different, but close parameters, $(c_1,c_2,c_3)=(-0.5, -0.7, -0.7)$ and $(c_1,c_2,c_3)=(-0.5, -0.7, -0.6)$ respectively, both of which are not in the multi-M set. The figures suggest a disconnected (although not totally disconnected) multi-J set in the first case, and a connected multi-J set in the second case. This implies that in this case the Fatou-Julia dichotomy fails in its traditional form -- and that the statement relating boundedness of the critical orbit with connectedness of the multi-J set does not hold for real networks. More precisely, we found parameters for which the multi-Julia set appears to be connected, although the critical multi-orbit is unbounded. On the other hand, the counterpart of the theorem may still hold, in the following form: ``The multi-J set is connected if the parameter belongs to the multi-M set.''

Part of our current work consists in optimizing the numerical algorithm for multi-M and J sets in real networks, with high enough resolution to allow (1) observation of possible fractal properties of multi-J sets and of multi-M sets boundaries and (2) computation of the genus of the filled multi-J sets, in attempt to phrase a topological extension of the theorem that takes into account the number of handles and tunnels that open up in these sets as their connectivity breaks down when leaving the Mandelbrot set.

\section{Discussion}
\label{discussion}

\subsection{Comments on our results}

In this paper, we used a combination of analytical and numerical approaches to propose possible extensions of Fatou-Julia theory to networked complex maps. We started by showing that, even in networks where all nodes are identical maps, their behavior may not be ``synchronized:'' the node-wise Mandelbrot sets may be identical in some cases, which in others they may differ substantially, depending on the coupling pattern. We then investigated how specific changes in the network hard-wiring trigger different effects on the structure of the network Mandelbrot and Julia sets, focusing in particular on observing topological properties (connectivity) and fractal behavior (Haussdorff dimension). We found instances in which small perturbations in the strength of one single connection may lead to dramatic topological changes in the asymptotic sets, and instances in which these sets are robust to much more significant changes.

More generally, our paper suggests a new direction of study, with potentially tractable, although complex mathematical questions. While existing results do not apply to networks in their traditional form, it appears that connectivity of the newly defined uni-Julia sets may still be determined by the behavior of the critical orbit. We conjectured a weaker extension of the Fatou-Julia theorem, which was based only on numerical inspection, and which remains subject to a rigorous study that would support or refute it.

\subsection{Future work}

There are a few interesting aspects which we aim to address in our future work on iterated networks. For example, we are interested in studying the structure of equi-M and uni-J sets for larger networks, and in understanding the connection between the network architecture and its asymptotic dynamics. This direction can lead to ties and applications to understanding functional networks that appear in the natural sciences, which are typically large. 

The authors' previous work has addressed some of these aspects in the context of continuous dynamics and coupled differential equations. However, when translating network architectural patters into network dynamics, the great difficulty arises from a combination of the graph complexity and the system's intractable dynamic richness. Addressing the question at the level of low-dimensional networks can help us more easily identify and pair specific structural patterns to their effects on dynamics, and thus better understand this dependence. The next natural step is to return to the search for a similar classification in high-dimensional networks, where specific graph measures or patters (e.g. node-degree distribution, graph Laplacian, presence of strong components, cycles or motifs) may help us, independently or in combination, classify the network's dynamic behavior.

\begin{figure}[h!]
\begin{center}
\begin{minipage}{0.26\textwidth}
\includegraphics[width=\textwidth]{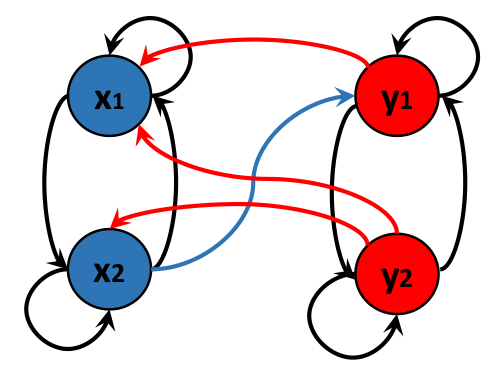}
\end{minipage}
\quad
\begin{minipage}{0.4\textwidth}
$\ds A$ = $\left[ \begin{array}{c|c}  M & A_1\\  \hline A_2 & M \end{array} \right]$, $\ds M = \left( \begin{array}{cc}  1 & 1\\ 1 & 1 \end{array} \right)$\\\\
$\ds A_1$ = $\left( \begin{array}{cc}  0 & 0\\  1 & 0 \end{array} \right)$, $\ds A_2 = \left( \begin{array}{cc}  1 & 0\\ 1 & 1 \end{array} \right)$
\end{minipage}
\begin{minipage}{0.27\textwidth}
\includegraphics[width=0.9\textwidth]{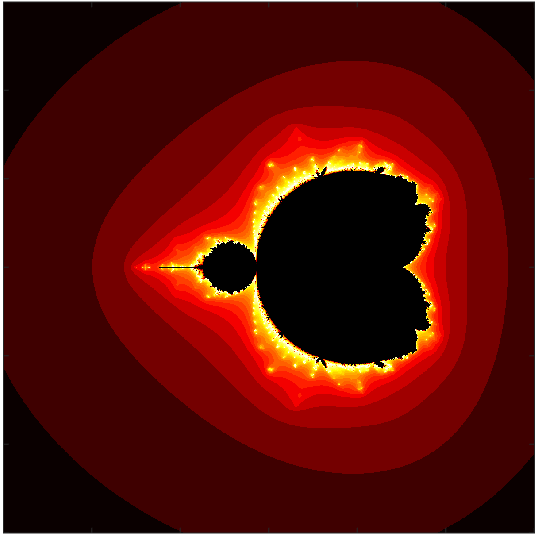}
\end{minipage}
\\
\begin{minipage}{0.26\textwidth}
\includegraphics[width=\textwidth]{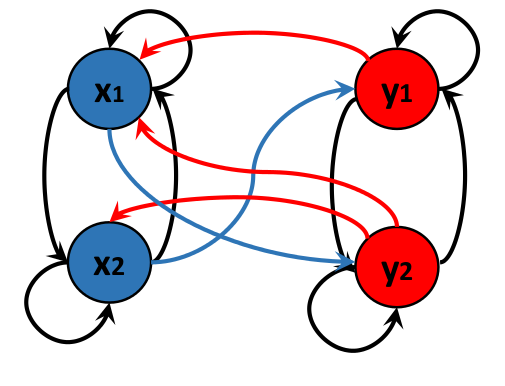}
\end{minipage}
\quad
\begin{minipage}{0.4\textwidth}
$\ds A$ = $\left[ \begin{array}{c|c}  M & A_1\\  \hline A_2 & M \end{array} \right]$, $\ds M = \left( \begin{array}{cc}  1 & 1\\ 1 & 1 \end{array} \right)$\\\\
$\ds A_1$ = $\left( \begin{array}{cc}  0 & 1\\  1 & 0 \end{array} \right)$, $\ds A_2 = \left( \begin{array}{cc}  1 & 0\\ 1 & 1 \end{array} \right)$
\end{minipage}
\begin{minipage}{0.27\textwidth}
\includegraphics[width=0.9\textwidth]{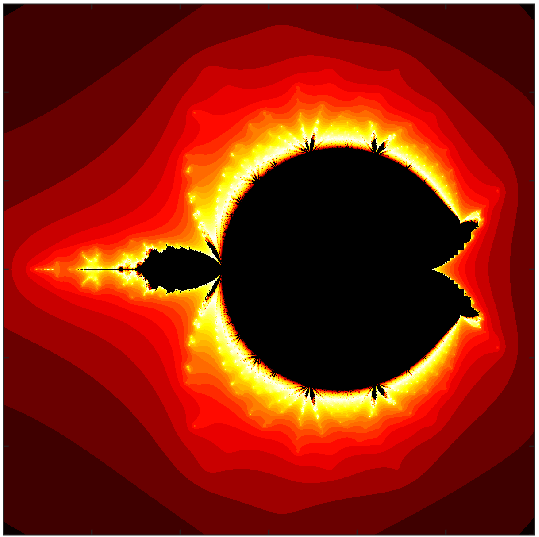}
\end{minipage}
\end{center}
\caption{\small \emph{{\bf Equi-M sets for two networks of size $N=4$} described schematically on the left, together with their adjacency matrices . Both systems have connectivity parameters $g_{xx} = g_{yy}=1/2$, $g_{xy} = g_{yx}=-1/2$.}}
\label{4dim_Mand}
\end{figure}

\begin{figure}[h!]
\begin{center}
\includegraphics[width=0.24\textwidth]{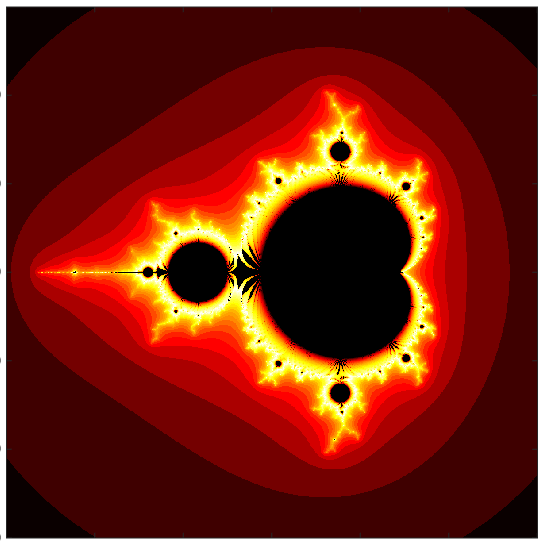}
\includegraphics[width=0.24\textwidth]{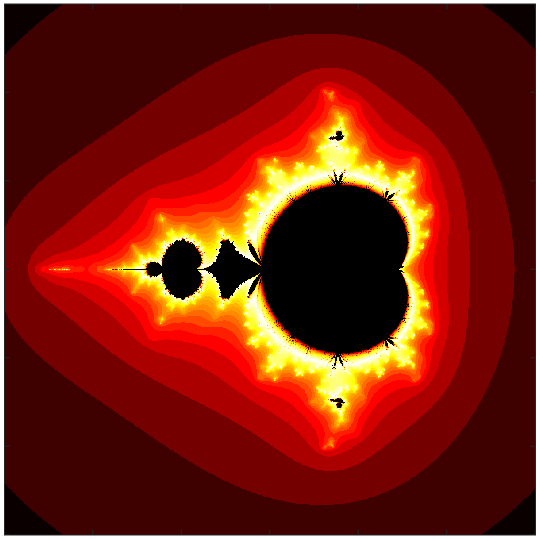}
\includegraphics[width=0.24\textwidth]{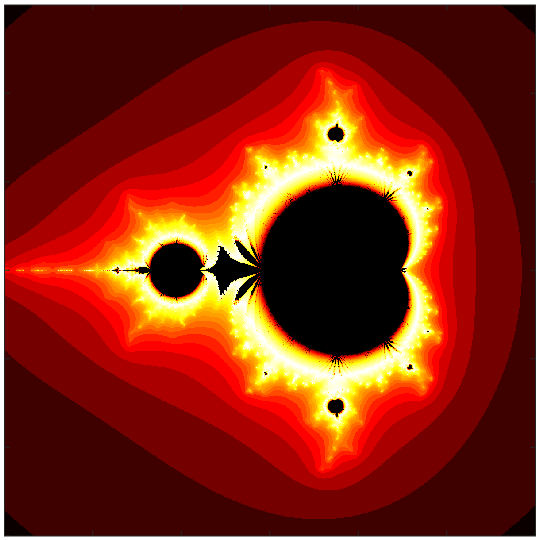}
\includegraphics[width=0.24\textwidth]{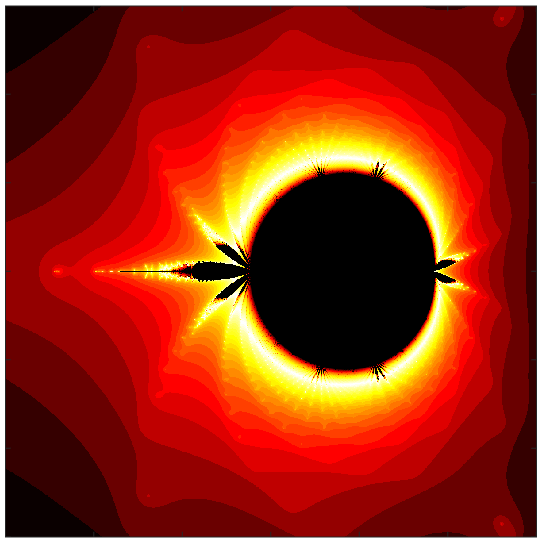}
\end{center}
\caption{\small \emph{{\bf Equi-M sets for a bipartite network with $20$ nodes}, formed of two cliques $X$ and $Y$, with $N=10$ nodes in each. Th adjacency matrix is therefore similar to those in Figure~\ref{4dim_Mand}, with square blocks $M$, $A_1$ and $A-2$ of size $N=10$. The densities (number of ones in each block, i.e. number of $X$-to-$Y$ and respectively $Y$-to-$X$ connecting edges) were takes in each panel to be (out of the total of $N^2=100$: {\bf A.} $N_{xy} = N_{yx}=10$; {\bf B.} $N_{xy}=15$, $N_{yx}=10$; {\bf C.} $N_{xy} = N_{yx}=15$; {\bf D.} $N_{xy} = N_{yx}=50$. In all cases, the connectivity parameters (i.e., edge weights) were $g_{xx}=g_{yy}=1/10$ and $g_{xy}=g_{yx}=-1/10$.}}
\label{10dim_Mand}
\end{figure}

Of high interest are methods that can identify robust versus vulnerable features of the graph from a dynamics standpoint. As Figures~\ref{4dim_Mand} and~\ref{10dim_Mand} show, it is clear that a small perturbation of the graph (e.g., adding a single edge) have the potential, even for higher dimensional networks, to produce dramatic changes in the asymptotic dynamics of the network, and readily lead to substantially different M and J sets. However, this is not consistently true. We would like to understand whether a network may have a priori knowledge of which structural changes are likely to produce large dynamic effects. This is a real possibility in large natural learning networks, including the brain -- where such knowledge probably affects decisions of synaptic restructuring and temporal evolution of the connectivity profile.

Our future work includes understanding and interpreting the importance of this type of results in the context of networks from natural sciences. One potential view, proposed by the authors in their previous joint work, is to interpret iterated orbits as describing the temporal evolution of an evolving system (e.g., copying and proofreading DNA sequences, or learning in a neural network). Along these lines, an initial $z_0$ which escapes to $\infty$ under iterations may represent a feature of the system which becomes in time unsustainable, while an initial $z_0$ which is attracted to a simple periodic orbit may represent a
feature which is too simple to be relevant or efficient for the system. Then the points on the boundary between these two behaviors (i.e., the Julia set) may be viewed as the optimal features, allowing the system to perform its complex function. We study how this ``optimal set of features'' changes when perturbing its architecture. 

Once we gain enough knowledge of networked maps for fixed nodes and edges, and we formulate which applications this framework may be appropriate to address symbolically, we will allow the nodes' dynamics, as well as the edge weights and distribution, to evolve in time together with the iterations. This process may account for phenomena such as learning, or adaptation --  a crucial aspect that needs to be understood about systems. This represents a natural direction in which to extend existing work by the authors on random iterations in the one-dimensional case.

\section*{Acknowledgements}

The work on this project was supported by the SUNY New Paltz Research Scholarship and Creative Activities program. We additionally want to thank Sergio Verduzco-Flores, for his programing suggestions, and Mark Comerford, for the useful mathematical discussions.

\bibliographystyle{plain}

\clearpage
\section*{Appendix A: Uni-J sets for higher dimensional networks}

\noindent The figures show four uni-J sets, for $N=4$ and $N=8$ nodes. The equi-parameters, adjacency matrices, and connectivity parameters of each network are given below, from left to right:\\

\begin{center}
\includegraphics[width=\textwidth]{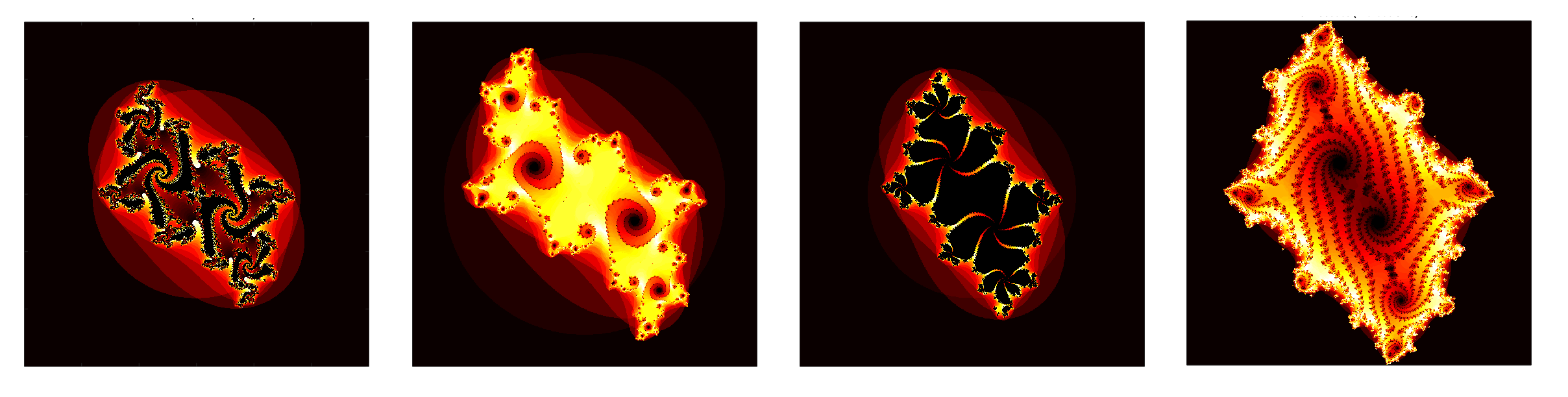}\\
\end{center}

\noindent {\bf Figure 1.} Size $N=4$, equi-parameter $c = -0.117-0.76i$.

\noindent Adjancy: $\ds A$ = $\left[ \begin{array}{c|c}  M & A_1\\  \hline A_2 & M \end{array} \right]$,
where $\ds M = \left( \begin{array}{cc}  1 & 1\\ 1 & 1 \end{array} \right)$, $\ds A_1$ = $\left( \begin{array}{cc}  1 & 0\\  0 & 0 \end{array} \right)$, $\ds A_2 = \left( \begin{array}{cc}  1 & 0\\ 1 & 0 \end{array} \right)$

\noindent Connectivity parameters: $g_{xx} = g_{yy}=1/2$, $g_{xy} = g_{yx}=-1/2$\\ \\

\noindent {\bf Figure 2.} Size $N=4$, equi-parameter $c = -0.117-0.856i$.

\noindent $\ds A$ = $\left[ \begin{array}{c|c}  M & A_1\\  \hline A_2 & M \end{array} \right]$,
where $\ds M = \left( \begin{array}{cc}  1 & 1\\ 1 & 1 \end{array} \right)$, $\ds A_1$ = $\left( \begin{array}{cc}  0 & 1\\  0 & 0 \end{array} \right)$, $\ds A_2 = \left( \begin{array}{cc}  1 & 0\\ 1 & 1 \end{array} \right)$

\noindent Connectivity parameters: $g_{xx} = g_{yy}=1/2$, $g_{xy} = g_{yx}=-1/2$\\ \\

\noindent {\bf Figure 3.} Size $N=4$, equi-parameter $c=-0.5622-0.62i$.

\noindent $\ds A$ = $\left[ \begin{array}{c|c}  M & A_1\\  \hline A_2 & M \end{array} \right]$,
where $\ds M = \left( \begin{array}{cc}  1 & 1\\ 1 & 1 \end{array} \right)$, $\ds A_1$ = $\left( \begin{array}{cc}  0 & 1\\  0 & 0 \end{array} \right)$, $\ds A_2 = \left( \begin{array}{cc}  1 & 0\\ 1 & 0 \end{array} \right)$

\noindent Connectivity parameters: $g_{xx} = g_{yy}=1/2$, $g_{xy} = g_{yx}=-1/2$\\ \\

\noindent {\bf Figure 4.} Size $N=8$, equi-parameter $c=-0.62-0.62i$.

\noindent $\ds A$ = $\left[ \begin{array}{c|c}  M & A_1\\  \hline A_2 & M \end{array} \right]$,
where $\ds M = \left( \begin{array}{cccc}  1 & 1 & 1 & 1\\ 1 & 1 & 1 & 1 \\ 1 & 1 & 1 & 1 \\ 1 & 1 & 1 & 1 \end{array} \right)$,  $\ds A_1 = \left( \begin{array}{cccc}  1 & 0 & 1 & 1\\ 1 & 1 & 1 & 0 \\ 1 & 0 & 0 & 0 \\ 0 & 1 & 0 & 0 \end{array} \right)$, $\ds A_2 = \left( \begin{array}{cccc}  0 & 0 & 0 & 0\\ 0 & 1 & 1 & 1 \\ 1 & 1 & 0 & 1 \\ 0 & 1 & 1 & 0 \end{array} \right)$

\noindent Connectivity parameters: $g_{xx} = g_{yy}=1/4$, $g_{xy} = g_{yx}=-1/4$\\ \\

\end{document}